\documentclass[12pt]{amsart}

\def\@typesizes{%
       \or{5}{6.5}\or{6}{7.5}\or{7}{8.5}\or{8}{11}\or{9}{12}%
       \or{10}{13}
       \or{\@xipt}{14}\or{\@xiipt}{15}\or{\@xivpt}{18}%
       \or{\@xviipt}{20}\or{\@xxpt}{24}}

\usepackage{amsthm}
\usepackage{amsmath}
\usepackage{amssymb}
\usepackage{graphicx}
\usepackage{enumerate}
\usepackage{color}
\allowdisplaybreaks
\numberwithin{equation}{section}
\numberwithin{figure}{section}

\theoremstyle{plain}
\newtheorem{theorem}{ Theorem}[section]
\newtheorem{proposition}[theorem]{ Proposition}
\newtheorem{lemma}[theorem]{ Lemma}
\newtheorem{corollary}[theorem]{ Corollary}
\newtheorem{example}[theorem]{ Example}
\newtheorem{remark}[theorem]{ Remark}
\newtheorem{definition}[theorem]{ Definition}
\newtheorem{conjecture}{ Conjecture}

\numberwithin{equation}{section}

\def\BET{\begin{theorem}}
\def\ENT{\end{theorem}}
\def\BEP{\begin{proposition}}
\def\ENP{\end{proposition}}
\def\BEL{\begin{lemma}}
\def\ENL{\end{lemma}}
\def\BEC{\begin{corollary}}
\def\ENC{\end{corollary}}
\def\BEE{\begin{example} \rm}
\def\ENE{\end{example}}
\def\BER{\begin{remark} \rm}
\def\ENR{\end{remark}}
\def\BED{\begin{definition} \rm}
\def\END{\end{definition}}
\def\BECJ{\begin{conjecture}}
\def\ENCJ{\end{conjecture}}

\def\bea{\begin{eqnarray}}
\def\eea{\end{eqnarray}}

\def\beas{\begin{eqnarray*}}
\def\eeas{\end{eqnarray*}}

\def\beq{\begin{equation}}
\def\eeq{\end{equation}}

\def\beal{\begin{align*}}

\def\eeal{ \end{align*} }

\def\row{ \nonumber \\ & & }
\def\roweq{\nonumber \\ &=& }
\def\rowleq{\nonumber \\  & \leq & }

\def\bfC{{\bf C}}

\def\bbA{{\mathbb A}}

\def\bbC{{\mathbb C}}
\def\bbD{{\mathbb D}}

\def\bbR{{\mathbb R}}

\def\cA{{\mathcal A}}

\def\cD{{\mathcal D}}

\def\cW{{\mathcal W}}

\newcommand{\WA}{{\cW}_{\cA}}
\newcommand{\wD}{{\widehat \cD}}

\def\ef{\eqref}
\def\wh{\widehat}
\def\wo{{\wh \omega}}

\begin{document}

\title{Quasiconformal symbols and  projected composition operators}

\author[S\"{o}nmez]{Sinem Yelda S\"{o}nmez}

\address{%
Faculty of Engineering and Architecture\\
Altınba\c{s} University\\
Ba\c{g}cilar, Istanbul\\
34218,Turkey}

\email{sinem.sonmez@altinbas.edu.tr} 
\address{Department of Mathematics and Statistics \\
University of Helsinki \\
P.O.Box 68 (Pietari Kalmin katu 5) \\
00014 Helsinki, Finland}
\email{sinem.sonmez@helsinki.fi}

\author[Taskinen]{Jari Taskinen} 
\address{Department of Mathematics and Statistics \\
University of Helsinki \\
P.O.Box 68 (Pietari Kalmin katu 5) \\
00014 Helsinki, Finland}
\email{jari.taskinen@helsinki.fi}




\begin{abstract}
We study projected composition operators $K_\varphi$ with quasiconformal symbols $\varphi$ on weighted Bergman spaces on the open unit disc
$\bbD$. If the symbol were conformal, i.e. a M\"obius transform of $\bbD$, the corresponding composition operator
would be automatically invertible at least in standard weighted spaces. We show that the invertibility remains, if the
Beltrami coefficient is small enough, in particular, it satisfies a certain vanishing condition  at the boundary
of the disc. 
We also consider the invertibility of $K_\varphi$ for symbols $\varphi$ which are conformal in an annulus 
$\{ R < |z| < 1 \}$. The weight classes in our considerations include both standard and exponentially decreasing weights. 
\end{abstract}

\maketitle

\section{Introduction}
\label{sec1}

The set of conformal mappings $\varphi$ from the the unit disk $\bbD$ of the complex plane $\bbC$ onto itself is very 
restricted, as it only consists of M\"obius-transforms. Among many consequences of this fact we remark that it often makes the 
study of the invertibility of analytic composition operators in spaces of analytic functions on $\bbD$ quite uninteresting, 
although  the question of the invertibility is in general most important for a given linear operator. However, a long time 
ago it has been observed that, in many other questions and phenomena in mathematical analysis, the analyticity of $\varphi$ 
can be relaxed and especially replaced by quasiconformality, yet the mapping still has many important properties similar 
to those of conformal mappings; see \cite{A}, \cite{Kos}, \cite{Vai}.

In the Bergman space $A^2$ of square integrable analytic functions $f: \bbD \to \bbC$, a composition operator $C_\varphi : f \mapsto f \circ \varphi$ with a quasiconformal symbol would, in general, map $f$ outside the space,
but this can be adjusted most naturally by composing $C_\varphi$ with the Bergman projection $P:L^2 \to A^2$.
This leads to the definition of the projected composition operator $K_\varphi = P C_\varphi$. In view of the
fact that the set of quasiconformal self-maps of the disc $\bbD$ is much richer than the set of 
conformal mappings (which induce invertible composition operators), we are naturally led to study the invertibility of 
operators $K_\varphi$ with quasiconformal $\varphi$ and, in particular, to derive sufficient conditions for the 
invertibility. This is the purpose of the present paper. 

In Theorem \ref{th3.1} we will prove the invertibility of projected composition operators $K_\varphi$ for
quasiconformal $\varphi$ with small enough Beltrami coefficients in weighted Bergman $L^2$-spaces. 
A function $\varphi:\bbD \to \bbD$ is quasiconformal, if it is an orientation-preserving homeomorphism, which belongs to 
the Sobolev space $W^{1,2}$ of functions whose first-order distributional derivatives belong to $L^2$ and if it is a 
weak solution of the Beltrami equation
\bea
\bar \partial f(z) = \mu(z) \partial f(z), \ \ \ z \in \bbD,  \label{1.2}
\eea
for some Lebesgue measurable function $\mu$ with $\sup_{z \in \bbD} |\mu(z)| < 1$, called the Beltrami coefficient. 
Here, $\partial = \partial / \partial  z$ and $\bar \partial = \partial / \partial \bar z$ are the Wirtinger derivatives; see the end of this section  for related terminology.  Apparently, all conformal mappings satisfy \ef{1.2}, as they are
solutions of the homogeneous $\bar \partial$-equation. 
The magnitude of the Beltrami coefficient gives a measure for the distance of a quasiconformal $\varphi$ to analyticity.
Accordingly, the assumptions of Theorem \ref{th3.1} include smallness requirements for $\mu$, see \ef{3.2}--\ef{3.3a}.
Another invertibility result for operators $K_\varphi$ in weighted $L^p$-spaces is contained in Proposition \ref{th4.1}. Here, 
the symbol  $\varphi: \bbD \to \bbD$ is assumed to be conformal in an annulus $\bbA_R= \{ R < |z| < 1 \}$ for a large enough
$R \in (0,1)$, and the proof is based on a different perturbation argument.

We refer to \cite{Kos}  for an introduction to the topic of quasiconformal mappings.

We will study the above mentioned questions in the setting of Bergman spaces with weighted $L^p$-norms, 
$1 < p< \infty$. Here, a weight $\omega : \bbD \to (0,\infty)$ is a continuous function which is radial 
($\omega(z)=\omega(\vert z\vert)$ for $z\in \mathbb{D}$) and satisfies $\lim_{r \to 1^-} \omega(r)= 0$. The classes of weights under consideration contain standard  and exponentially decreasing weights, and they will be discussed in
detail in the next section. Given $1 <p<\infty$ and a weight $\omega$, we denote by $L_\omega^p$ the $L^p$-space with respect to 
the measure $\omega dA$,  where $dA$ is the normalized Lebesgue area measure  on $\bbD$, i.e. $dA = \pi^{-1} r dr d\theta$, and write 
\begin{equation*}
    \|f\|_{p,\omega}^p := \int\limits_{\mathbb{D}} \vert f(z)\vert^p \omega(z) dA(z) .
\end{equation*}
The weighted Bergman space $A_\omega^p$ is the closed subspace of $L_\omega^p$ consisting of analytic functions. Our 
considerations also contain the case with no weight, i.e. $\omega =1$ on $\bbD$, in which case we denote the 
spaces and the norm by $L^p, A^p$ and $\Vert \cdot \Vert_p$, respectively. 

As for linear operators in these Banach spaces, we denote by $P_\omega$ the orthogonal projection from $L_\omega^2$ onto 
$A_\omega^2$. In the sequel, we will consider situations with two weights $\omega$ and $\nu$, which may or may not be the 
same, such that $P_\omega$ is a bounded projection from $L_\nu^p$ onto $A_\nu^p$. 
Given a map $\varphi: \bbD \to \bbD$ and a weight $\omega$, we define the projected composition operator 
\begin{equation*}
K_{\omega,\varphi} f=P_\omega C_\varphi = P_\omega (f\circ\varphi), 
\end{equation*} 
where $C_\varphi : f \mapsto f \circ \varphi$ is a composition operator. We also write
$K_\varphi = K_{\omega,\varphi}$, if the weight is clear from the context. 

Projected composition operators were probably first studied in the setting of analytic function spaces in the paper 
\cite{Mishra} by S. Pattanyak, C. Mohapatra and A. Mishra. Their boundedness properties  have been considered for
example by R. Rochberg \cite{Rochberg} in Hardy spaces and by C. Zhao \cite{Zhao} and 
$\check{\mbox{Z}}$. $\check{\mbox{C}}$u$\check{\mbox{c}}$kovi\'c \cite{Cuckovic} 
in Bergman spaces. Also, the reader is probably aware of the fact that there exists a vast literature
on the properties of analytic composition operators, that is, on operators with analytic $\varphi$. A review of 
that field of research is outside the scope of this article. Finally, it should be mentioned that 
composition operators with quasiconformal symbols acting in Sobolev spaces have been investigated in many 
papers, see for example \cite{Ast},  \cite{Kos14}, \cite{Rei}. However, due to the  different nature of Sobolev 
and Bergman  spaces, the interesting questions and methods in the mentioned papers are quite different from the 
present study. The recent paper \cite{Fang} treats the boundedness of the operator $C_\varphi$ from the Bergman
space into $L^2$, however, only in the unweighted case.

As for the contents of his paper,  the main results on the invertibility of projected composition 
operators are presented in Theorem \ref{th3.1} and also Proposition \ref{th4.1}.
Section \ref{sec2} also contains the preliminary material like the definitions of the weight classes,  a simple 
construction of a left inverse for the $\bar \partial$-operator in a subspace of 
$L_\omega^p$ (Proposition \ref{prop 2.4}), and
some remarks on the boundedness of composition operators with non-analytic symbols. In Section \ref{sec5} we
present a related example of a non-invertible projected composition operator with a quasiconformal symbol. 

We will use standard notation, see for example \cite{Zhu}. All function spaces are defined on the domain
$\bbD$, unless otherwise indicated. 
Moreover, $\Delta$ denotes the Laplace operator in the real  
coordinates of the plane. If $z = x + iy \in \bbC$, the Wirtinger derivatives are
defined as $\partial = \partial / \partial z = ( \partial / \partial x - i \partial / \partial y)/2$
and $\bar \partial = \partial / \partial z = ( \partial / \partial x + i \partial / \partial y)/2$.
If $\varphi \in W^{1,2}$  is a  homeomorphism of $\bbD$ onto itself, we denote by $J_\varphi$ its Jacobian
determinant (in the real planar coordinates).

If $f$ and $g$ are real valued, posit\-ive expressions depending for example on a
parameter $x$, the expression $f(x) \asymp g(x)$ means that  there is constant $C > 0$ such that
$C^{-1} f(x) \leq g(x) \leq Cf(x)$ for all $x$.

\section{Preliminaries on weight classes and the $\bar \partial$-equation in weighted $L^p$-spaces}
\label{sec2}

We present in this section the definitions of the weight classes  and some basic results related with them and then 
proceed to a simple construction of a left inverse of the $\bar \partial$-operator in weighted $L^p$-spaces. 
Our weight classes contain both standard ($1^\circ$, below) and exponentially decreasing $(3^\circ$)
weights. Also, upper doubling weights ($2^\circ$, below) will be mentioned in some relevant preliminary results, although
we will not formulate the main invertibility results for them. 

As mentioned, by a weight in general we mean a continuous, 
radial function $\omega: \bbD \to (0,\infty)$ with limit 0 at the boundary of $\bbD$. 

\BED \label{def1.1}
\noindent 
$1^\circ$. A weight $\omega$ is standard, if  $\omega(z) = (1-|z|^2)^\alpha$ 
for some $\alpha \geq 0$. In particular, the unweighted or constant weight case corresponds to $\alpha = 0$.

\smallskip

\noindent
$2^\circ$. A weight $\omega$ is upper doubling and belongs to the class 
$\wD$, if the weight 
$\widehat \omega (z) = \int_{|z|}^1 \omega(s) ds $ satisfies $\wo (r) \leq C \wo 
(\frac{1+r}{2})$, for some constant $C > 0$ and for all $r \in (0,1)$. 
See  \cite{AS}, Section 4 and \cite{PRAdv}. Note that all standard weights  belong to $\wD$.

\smallskip

\noindent
$3^\circ$. A weight $\omega$ belongs to the class $\WA$, if it is  of the form
\bea
\omega (z) = e^{-2 \phi(z)} , \ \ \ z \in \bbD,  \label{2.0}
\eea
for some subharmonic, twice continuously differentiable function $\phi$, which 
satisfies $\left(\Delta\phi(z)\right)^{-1/2}\asymp \tau(z)$ for some continuously differentiable, 
radial and  positive function $\tau$ with the following properties: 

\noindent $(i)$ $\lim_{r \to 1^- } \tau(r) =0$ and $\lim_{r \to 1^- } \tau'(r) =0$;  

\noindent $(ii)$ there is 
a constant $C>0$ such that $\tau(r)(1-r)^{-C}$ increases for $r$ close to $1$, or
$\tau'(r)\log \left(\tau(r)\right) \rightarrow 0$ as $r\to 1$; 

\noindent $(iii)$ there holds
\bea
\varphi(r) = \frac{\partial \varphi(r) }{\partial r} \asymp \frac1{R(r)} , \ \ r\in(0,1),
\eea
where
\begin{equation}
 R(r)= \frac{\tau(r)^2}{1-r} , \  \ \ \ r \in (0,1).   \label{def:R}
\end{equation} 
with a radial extension $R(z):= R(|z|)$ for $z \in \bbD$. (Note that 
$R$ tends to 0 at the boundary of $\bbD$.)
\END


There holds $\wD \cap \WA = \emptyset$. 
The weight class $\WA$ was introduced in \cite{ATV} in the study of pointwise estimates for the Bergman kernel
in large Bergman spaces. The prototype of the weights in $\WA$ is the exponential type weight
\begin{equation}\label{eqn:Exp}
	\omega_{a,b}(z) = \exp\left(\frac{-b}{(1-|z|^2)^a}\right),   
	\quad  b,a > 0 .
\end{equation} 
The class $\WA$ is contained in a slightly larger class $\mathcal{W}$, which together with other
related classes of rapidly decreasing weights appear in the study of various  operators on $A^p_\omega$, such as 
Bergman projections \cite{HLS}, Hankel \cite{GP1, HP2021}, integration \cite{const2010, D2, PP1, PP2}
and  Toeplitz operators \cite{APP, ST} and others.

Given a weight $\omega$ and $1 < p < \infty$, we define the  weight 
\begin{equation}
\nu_p = \left\{ 
\begin{array}{ll}
\omega, \ \ &\mbox{if} \ \omega \in \wD , \\
\omega^{p/2} , \ \ &\mbox{if} \ \omega \in \WA. 
\end{array}
\right.  \label{1.8}
\end{equation}
With a small  abuse of notation,  we  denote by $L_\nu^p$ and $A^p_\nu$ the $L^p$- and Bergman spaces 
with respect to the measure  $\nu_p \,dA$. 
The need of such a notation is explained in the following result, which clarifies the question on the boundedness of 
Bergman projections in  relation with our weight classes. If $\omega \in \wD$, the result is proved in \cite{PRAdv}, Theorem 7, 
and if $\omega \in \WA$, it follows from Theorem 4.1 of \cite{HLS}, since our class $\WA$ is contained
in the class $\cW_0$ of the citation. 

\BEL  \label{lem2.0}
Assume $\omega \in \wD \cup \WA$ and  $1 < p < \infty$. Then, the Bergman projection
$P_\omega$ is a bounded operator from $L_\nu^p$ onto $A_\nu^p$ for all $p \in (1,\infty)$. 
\ENL

We denote by $d_P = d_P(p,\omega)$ the norm of $P_\omega$ in $L_\nu^p$. 

Let us next recall a Littlewood-Paley-type inequality, which will be needed later. 

\BEL
\label{lem2.1}
Assume  $\omega$ is a standard weight or belongs to $\WA$, and let us define for all $z \in\bbD$
\begin{equation}
\rho_2 (z) = \left\{ 
\begin{array}{ll}
(1- |z|^2)^2, \ \ &\mbox{if} \ \omega \ \mbox{is standard}  , \\
R^{2}(|z|) , \ \ &\mbox{if} \ \omega \in \WA.
\end{array}
\right.  \label{1.8a}
\end{equation}
There exists a constant $d_{LP}= d_{LP}(p,\omega) > 0$ such that the inequality
\bea
\int\limits_\bbD | f'|^2 \rho_2 \omega \, dA \leq d_{LP}^2 \Vert f \Vert_{2, \nu}^2  
\label{2.1}
\eea
holds for all $f \in A_\omega^2$. 
\ENL

For a standard $\omega $ this claim is well-known, see for example the remark after the proof of Theorem 4.28 in
\cite{Zhu}; see also \cite{PRAdv}, Theorem 5, for a more general result.  
If $\omega \in\WA$, \ef{2.1} follows from Theorem 5.1.$(i)$ of \cite{ATV}: one takes $p=q=2$, $n=1$ and the 
unweighted Lebesgue measure $dA$ for the Borel measure $\mu$. Then, in the notation of \cite{ATV}, 
$\mu(D_\delta)(z) \asymp \tau^2(z)$ for $z \in \bbD$, and the claim follows from the citation.

\bigskip

Given $p, \omega$ and $\nu_p$ as in \ef{1.8}, we define the subspace  $H \subset L_\nu^p$ to 
consist of functions $g$ such that  the expression 
\bea
\Vert g \Vert_H^p = \int\limits_\bbD |g|^p \nu_p \, dA +  \int\limits_\bbD |\bar\partial g|^p
\nu_p \, dA   \label{2.6}
\eea
is finite, where $\bar \partial g$ denotes the weak derivative. Then, the space $H$ endowed with the
norm $\Vert \cdot \Vert_H$ is a Banach space. 

We will need the following simple observations. 

\BEL \label{lem2.2}
Let $1 < p \leq 2$ and assume that $\omega$ is a standard weight or belongs to the class $\WA$. The operator $\bar \partial$ is a bounded surjection $H \to L_\nu^p $.
\ENL

Proof.  Assume first $\omega \in \WA$. Then, the function $\phi$ in \ef{2.0} is  subharmonic, and the case $p=2$ in 
Lemma \ref{lem2.2} follows from H\"ormander's solution  of the $\bar \partial$ problem, see \cite{Ho}:
for every $f \in L_\nu^2 = L_\omega^2$ there exists $u \in L_\omega^2$ with $\bar \partial u = f$.
(Another reference is Theorem 4.2.1 in \cite{HoN}, where the weight factors $1 + |z|^2$ in the inequality $(4.2.8)'$ do 
not matter, since the domain $X = \bbD$ is bounded.) By definition, $u$ belongs to the space $H$.
If $1 < p < 2$, the function $\varphi/p$ is still subharmonic, and one can argue in the same way, by using the 
result of  \cite{FS} instead of \cite{Ho}, see also Theorem 1 in  \cite{BB}.  

If $\omega(z) = (1- |z|^2)^\alpha$ is a standard weight with $\alpha \geq 0$, the above argument still applies, since
$\omega = e^{- 2\phi}$ with the subharmonic function $\phi(z) = - \alpha \log(1-|z|^2)/2$. The unweighted case 
$\alpha =0$ is not excluded. 
\ \ $\Box$

\BEL  \label{lem2.3}
Let $p \in (1,2]$ and assume  $\omega$ is a standard weight or belongs to $\WA$. Then, the Bergman space 
$ A_\nu^p$ is a closed 
subspace of the Banach space $H$, and  $P_\omega$ is a bounded projection operator from $H$ onto $A_\nu^p$. 
The space $H$ has a direct sum decomposition $H = A_\nu^p \oplus Y$ with the closed subspace
(complement of $A_\nu^p$) $Y= Q_\omega (H)$,  where $Q_\omega := I - P_\omega$.
\ENL

Proof. Since $\bar \partial f = 0$ for all $f \in A_\nu^p$, we note that $A_\nu^p \subset H$ and
the norm $\Vert  \cdot \Vert_H$ coincides with the norm of $L_\nu^p$ in this subspace. We infer 
that $A_\nu^p $ is complete, when endowed with the norm of $H$, and thus a closed subspace. 
Moreover, the Bergman projection $P_\omega$ is a bounded operator from $H$ onto $A_\nu^p$, since
by Lemma \ref{lem2.0},  
\bea
\Vert P_\omega f \Vert_H  = \Vert P_\omega f \Vert_{p,\nu} \leq d_P \Vert  f \Vert_{p,\nu}
\leq d_P \Vert  f \Vert_H 
\eea 
for all $f \in H$. That $P_\omega$ is a projection onto $A_\nu^p$ follows trivially from its
mapping properties in $L_\nu^p$. 

The remaining statements belong to elementary functional analysis.\ \ $\Box$

\bigskip

The following result on the left inverse of the $\bar \partial$-operator is known, but we give a proof for the sake of 
the completeness of the presentation. Indeed, this one variable version is sufficient for our purposes, and its proof
is simpler than those for higher dimensional domains, see for example the review in Notes for Chapter V in \cite{Range}.

\BEP  \label{prop 2.4}
Assume  $p \in (1, 2]$ and that $\omega$ is a standard weight or belongs to $\WA$. 
There exists a bounded linear operator $M : L_\nu^p \to L_\nu^p $ 
such that $M \bar \partial f = Q_\omega f$ for all $f \in H$  and hence 
\bea
M  \bar \partial f = f \ \ \ \mbox{for all} \ f \in Y = Q_\omega(H) \subset H \subset  L_\nu^p . 
\label{2.10}
\eea
\ENP

We denote by $d_M = d_M(p,\omega)$ the operator norm of $M$ in the space $L_\omega^p$. 

\bigskip

Proof. We consider $\bar \partial$ as a bounded, surjective operator $H \to L_\nu^p$, see Lemma \ref{lem2.2}. There
holds ker\,$\bar \partial \cap Y = A_\nu^p \cap Y = \{ 0 \}\subset H$, hence, 
by the open mapping theorem, $\bar \partial$ is an isomorphism from $Y $ onto $ L_\nu^p$, and 
in particular it is bounded from below. Let us denote $M = (\bar 
\partial)^{-1} : L_\nu^2 \to Y \subset H$. Then, the operator $M \bar \partial$ is the identity 
on $Y$ and its kernel coincides with $A_\nu^p$. Thus, $M \bar \partial$ is the projection 
onto $Y$ corresponding to the direct sum decomposition $H = A_\nu^p\oplus Y$, i.e.
$ M \bar \partial = Q_\omega$.

Note that since $M$ is bounded as an operator $L_\nu^p \to Y \subset H$, it is  
also bounded $L_\nu^p \to L_\nu^p$. \ \ $\Box$ 

\bigskip
At the end of this section we briefly review the boundedness properties of composition operators;
boundedness of $C_\varphi$  will be needed in the proof of the invertibility results. Due to the non-analyticity of the symbols and the generality
of our weight classes, known  results on the boundedness of analytic composition operators do not suffice for our purposes. 
A sufficient condition for the boundedness of $C_\varphi$ is given in the following observation. Recall that $J_\varphi$
denotes the Jacobian.

\BEP \label{prop2.0}
Let $p \in (1,\infty)$ and $\omega \in \wD \cup \WA$. Assume $\varphi : \bbD \to \bbD$ is a homeomorphism such that $\varphi \in 
W^{1,2}$ and  $\varphi^{-1} =: \psi\in W^{1,2}$ and such that, for some constants $b_2 \geq b_1 >0$, 
\bea
b_1 \omega \circ \varphi (z) |J_\varphi(z)| \leq \omega(z) \leq b_2 \omega \circ \varphi (z) |J_\varphi(z)| \label{2.12}
\eea
for all $z \in \bbD$. Then, the composition operators $C_\varphi$ and $C_\psi$ are bounded as mappings
from $L_\omega^p$ into $L_\omega^p$. 
\ENP

For the operator $C_\varphi$, this follows by a simple change of variable, which is justified for example by  Lemma 5.12 of \cite{Kos}:
\beas
&& \Vert f \circ \varphi \Vert_{p,\omega }^p  =  \int\limits_\bbD |f \circ\varphi|^p \omega   dA 
\leq b_2  \int\limits_\bbD |f  \circ\varphi |^p \omega\circ \varphi  |J_\varphi |   dA 
= b_2  \int\limits_\bbD |f|^p \omega  dA   . 
\eeas
The operator $C_\psi$ can be treated in the same way, since the second inequality in \ef{2.12} for $\psi$ follows from the 
first inequality of \ef{2.12} for $\varphi$. Also, cf. \cite{GS}, Section 2. 

We finally mention the following characterization of the boundedness of $C_\varphi$, when considered as a 
mapping from the weighted Bergman space into the corresponding $L^p$-space.  First, 
if $1 < p<  \infty$ and $\omega \in \wD \cup \WA$ are given, we recall that a positive measure $dm$
on $\bbD$ is called a Carleson measure for the space $A_\omega^p$, if there exists a constant $C>0$
such that
\beas
\int\limits_\bbD |f|^p dm \leq C \int\limits_\bbD |f|^p \omega dA
\eeas
for all $f \in A_\omega^p$.  Moreover, if $\varphi:\bbD \to \bbD$ is Lebesgue measurable, we define the 
corresponding pull-back measure $V_\varphi = V_\varphi(p,\omega)$ by setting 
\beas
V_\varphi(E) = \int\limits_{\varphi^{-1}(E)} \omega dA 
\eeas
for all Lebesgue measurable sets $E \subset \bbD$.  We leave it to the reader to verify that the proof of 
Proposition 1 of \cite{Cuckovic} can be generalized to show the next statement.

\BEP  \label{prop2.5} Assume that $p \in (1, \infty)$, $\omega \in \wD \cup \WA$ and that 
$\varphi:\bbD \to \bbD$ is Lebesgue measurable. Then, the composition operator $C_\varphi $ is bounded $A_\omega^p 
\to L_\omega^p$, if and only if the pull-back measure $V_\varphi(p,\omega)$ is Carleson for the space 
$A_\omega^p$.
\ENP

%

\section{Invertibility of $K_\varphi$ in the Hilbert space case}  \label{sec3}

In this section we formulate and prove the main invertibility result for projected composition operators
with quasiconformal symbols $\varphi$. In the case of analytic composition operators, it is quite obvious that a 
M\"obius-transform $\varphi$ would generate an invertible operator with inverse $C_\psi$, where 
$\psi = \varphi^{-1}$, if the weight in the underlying space is standard or at least doubling. 
Theorem \ref{th3.1} generalizes this result for quasiconformal symbols with small enough Beltrami coefficients. 

Recall that if $p$ and $\omega $ are given, $d_M $  denotes  the operator norm of the operator
$M: L_\omega^2 \to L_\omega^2$, \ef{2.10}, and that the constant $d_{LP}$ was defined in \ef{2.1}.
In addition, for a homeomorphism $\varphi: \bbD \to \bbD$ with inverse $\psi = \varphi^{-1}$,
we denote by $d_\varphi = d_\varphi(p,\omega) $ and  $d_\psi = d_\psi(p,\omega)$ the operator norms of 
$C_\varphi : L_\omega^2 \to L_\omega^2$ and $C_\psi:  L_\omega^2  \to L_\omega^2$, respectively, if the operators 
are bounded.  
Also,   $\mu$  stands for  the Beltrami coefficient of a quasiconformal mapping $\varphi$, i.e., there holds 
\bea
\bar \partial \varphi = \mu  \partial \varphi .   \label{3.1}
\eea

\BET \label{th3.1}
Let  $\varphi: \bbD \to \bbD$ be a  quasiconformal mapping with inverse $\psi = \varphi^{-1}$ 
such that $|\mu(z)| \leq 1/2$ for all $z \in \bbD$ and such that the composition operators $C_\varphi$ and $C_\psi$ 
are bounded $L_\omega^2 
\to L_\omega^2$, and let the constant $\delta \in (0,1/\sqrt{2})$ be arbitrary. 

\noindent $(i)$ If $\omega$ is a standard weight $\omega(z) = (1- |z|^2)^\alpha$, $\alpha \geq 0$, and the 
Beltrami coefficient satisfies for all $z \in \bbD$, 
\bea
|\mu(z)| < \gamma_\psi \frac{(1- |\varphi(z)|^2)^{1+ \alpha/2}}{(1- |z|^2)^{\alpha/2}} \ \ \ \ 
\mbox{with the constant} \ \ 
\gamma_\psi = \frac{\delta}{d_{LP} d_M d_\psi},  \label{3.2}
\eea
and
\bea
|\mu(z)| < \gamma_\varphi \frac{(1- |\psi(z)|^2)^{2+ \alpha}}{(1- |z|^2)^{\alpha}} \ \ \ \ \mbox{with the constant} \ \ 
\gamma_\varphi = \frac{\delta}{d_{LP} d_M d_\varphi},  \label{3.2a}
\eea
then the projected composition operator $K_\varphi = K_{\omega,\varphi}: A_\omega^2 \to A_\omega^2 $ is invertible.

\noindent $(ii)$ If $\omega \in\WA $  and for all $z\in \bbD$ there holds
\bea
|\mu(z)| <  \gamma_\psi \frac{ \tau(\varphi(z))^2 \omega (\varphi(z))^{1/2}}{(1-|\varphi(z)|)\omega(z)^{1/2}} 
 \ \ \ \ \mbox{for all} \ z \in \bbD, \label{3.3}
\eea
and  
\bea
|\mu(z)| <  \gamma_\varphi \frac{ \tau(\psi(z))^2 \omega (\psi(z))^{1/2}}{(1-|\psi(z)|)\omega(z)^{1/2}}
 \ \ \ \ \mbox{for all} \ z \in \bbD, \label{3.3a}
\eea
then the projected composition operator $K_\varphi= K_{\omega,\varphi}:A_\omega^2 \to A_\omega^2 $ is invertible.
\ENT

\bigskip

Proof. We prove both cases $(i)$ and $(ii)$ simultaneously. Given $f \in A_\omega^2$, we start by showing that
$f \circ \varphi$ belongs to the space $H$ defined in \ef{2.6}. First, the assumptions of the theorem
yield
\bea
\Vert  f \circ \varphi  \Vert_{2,\omega}^2  \leq d_\varphi^2 \Vert f   \Vert_{2,\omega}^2  \label{J20}
\eea
To treat the derivative term in \ef{2.6}, we use the relation $\bar \partial (f \circ \varphi) = 
(\partial f \circ\varphi ) \bar \partial \varphi$ and \ef{3.1},  and calculate as follows: 
\bea
& & \Vert \bar \partial ( f \circ \varphi)  \Vert_{2,\omega}^2
= \int\limits_\bbD | \partial f \circ \varphi |^2 \, | \bar \partial \varphi|^2  \omega dA
\roweq 
\int\limits_{\bbD} | \partial f \circ \varphi |^2 \, \big( - 2 \mu^2 | \bar \partial \varphi|^2 
+ (1 + 2 \mu^2) | \bar \partial \varphi|^2  \big) \omega dA
\roweq 
\int\limits_{\bbD} | \partial f \circ \varphi |^2 \, \big( - 2 \mu^2 | \bar \partial \varphi|^2 
+ (1 + 2 \mu^2)  \mu^2  | \partial \varphi|^2  \big)  \omega dA
\roweq 
 \int\limits_{\bbD} | \partial f \circ \varphi |^2 \, \big( - 2 \mu^2  | \bar \partial \varphi|^2 
+ 2 \mu^2  | \partial \varphi|^2 + \mu^2  (2 \mu^2 -1) | \partial \varphi|^2  \big) \omega dA .
\label{3.4}
\eea
Recall that the Jacobian $J_\varphi$ of the function $\varphi: \bbD \to \bbD$, when understood as
a coordinate transform for real area integrals, equals 
\bea
J_\varphi = |\partial \varphi|^2 - |\bar\partial \varphi|^2 , \label{3.6}
\eea
see for example \cite{A}, formula (9) on p. 4, or \cite{Kos}, p. 80.  Moreover, it follows from \ef{def:R} and \ef{1.8a}  
that both \ef{3.2} and \ef{3.3} can be written as the inequality
\bea
\mu(z)^2 \omega(z) \leq \gamma_\psi^2 \rho_2(\varphi(z)) \omega(\varphi(z)) , \ \ z \in \bbD.
\eea
Hence, we can use \ef{3.6} and $2 \mu^2- 1 \leq 0$ and \ef{2.1} to infer that \eqref{3.4} is bounded by
\bea
& &  \int\limits_{\bbD} | \partial f \circ \varphi |^2 \, 2 \mu^2 \big( -  | \bar \partial \varphi|^2 
+  | \partial \varphi|^2  \big)\omega  dA 
\leq \int\limits_{\bbD} | \partial f \circ \varphi |^2 |J_\varphi|  2 \mu^2 \omega dA 
\rowleq
2 \gamma_\psi^2 \int\limits_{\bbD} | \partial f \circ  \varphi   |^2   \rho_2 \circ \varphi \, \omega \circ \varphi |J_\varphi| dA 
\leq
2\gamma_\psi^2 \int\limits_{\bbD} | \partial f  |^2  \rho_2 \,  \omega (z) dA  
\rowleq 
2\gamma_\psi^2 \int\limits_{\bbD} |f'  |^2   \rho_2 \,  \omega  dA  \leq 
2\gamma_\psi^2 d_{LP}^2  \Vert f \Vert_{2,\omega}^2   . \label{3.8}
\eea 
So, combining \ef{3.4}--\ef{3.8} yields
\bea
\Vert \bar \partial (f \circ \varphi) \Vert_{2,\omega}^2 
\leq 2\gamma_\psi^2 d_{LP}^2 \Vert f \Vert_{2,\omega}^2, \label{J22}
\eea
and together with \ef{J20} this  proves that $f \circ \varphi \in H$. Now, we can use Proposition \ref{prop 2.4}
to obtain 
\bea
& & \Vert C_\varphi f  - K_\varphi f\Vert_{2, \omega}^2 
= \Vert (P_\omega - I)C_\varphi  f  \Vert_{2,\omega}^2 
= \Vert Q_\omega  ( f \circ \varphi)  \Vert_{2,\omega}^2
\roweq 
\Vert M \bar \partial ( f \circ \varphi)  \Vert_{2,\omega}^2
\leq d_M^2 \Vert \bar \partial ( f \circ \varphi)  \Vert_{2,\omega}^2
\eea
This and  \ef{J22} yields us 
\bea
\Vert C_\varphi f  - K_\varphi f\Vert_{2, \omega}^2 \leq  2 \gamma_\psi^2 d_{LP}^2 d_M^2 
\Vert f \Vert_{2,\omega}. 
\label{3.10}
\eea
Since the inverse of a quasiconformal mapping is quasiconformal with the same Beltrami coefficient (\cite{Kos}, Theorem
6.3), we also obtain by using \ef{3.2a}, \ef{3.3a} instead of \ef{3.2}, \ef{3.3}, 
\bea
\Vert C_\psi f  - K_\psi f\Vert_{2, \omega}^2
= \Vert C_\psi f  - K_{\omega, \psi} f\Vert_{2, \omega}^2 \leq 2\gamma_\varphi^2 d_{LP}^2 d_M^2  \Vert f \Vert_{2,\omega} .  
\label{3.12}
\eea
for all $f \in A_\omega^2$, by the same calculation as in  \ef{J20}--\ef{3.10}.

Assume now that  $f \in A_\omega^2$ is arbitrary. We have 
$K_\psi C_\varphi f = P_\omega (f\circ \varphi \circ \psi) = P_\omega f = f $ due to the analyticity of $f$.
Hence, 
\bea
& & f - K_\psi K_\varphi f =  K_\psi (C_\varphi - K_\varphi) f . \label{3.14}
\eea
Here,  we estimate using \ef{3.10} and the choice of the constant $\gamma_\psi$ in \ef{3.12},
\beas
& & \Vert (I- K_\psi K_\varphi)f \Vert_{2,\omega} = \Vert K_\psi \big( C_\varphi - K_\varphi \big) f 
\Vert_{2, \omega}  
\leq d_\psi \Vert ( C_\varphi - K_\varphi ) f \Vert_{2, \omega}  
\leq \sqrt{2}\delta \Vert f \Vert_{2,\omega} . 
\eeas
which, in view of the choice $\delta < 1/\sqrt{2}$ and the Neumann series, 
means that the operator $K_\psi K_\varphi$ is invertible. In the same way,
interchanging the roles of $\varphi$ and $\psi$ in \ef{3.14} and using \ef{3.10}, 
one proves that  $K_\varphi K_\psi$ is invertible. We conclude that $K_\psi ( K_\varphi  K_\psi)^{-1}$ is the bounded 
inverse of $K_\varphi$:
\bea
& & K_\varphi \big(K_\psi (K_\varphi K_\psi)^{-1} \big) = I \ \ \ \ \mbox{and} \ \ \ \
\big(K_\psi (K_\varphi K_\psi)^{-1} \big) K_\varphi
\roweq
K_\psi (K_\varphi K_\psi)^{-1} K_\varphi K_\psi K_\varphi ( K_\psi K_\varphi)^{-1} 
= K_\psi K_\varphi ( K_\psi K_\varphi)^{-1} = I .  \ \ \Box
\eea

\bigskip

We refer to Proposition \ref{prop2.0} for a sufficient condition for the boundedness of the composition
operators, cf. the assumptions of the theorem. As it is obvious from the proof, it would be enough to assume that $C_\varphi$ 
is bounded as an
operator from the closed span of $A_\omega^2 \cup C_\psi (A_\omega^2)$ to $L_\omega^2$ and similarly for the operator $C_\psi$.

\section{On the invertibility of $K_\varphi$ in weighted $L^p$-spaces}  \label{sec4}

In this section we present another approach to the invertibility of
projected composition operators, which also works  in weighted Bergman spaces which are not necessarily
Hilbert. Instead of M\"obius mappings, the idea is to consider functions 
$\varphi : \bbD \to \bbD$ which only map the annulus  $\bbA_R = \bbD \setminus \bbD_R$, where $R \in (0,1)$ and 
$\bbD_R = \{z \in \bbD \, : \, |z| \leq R\}$, conformally onto a neighborhood of $\partial \bbD$. 

Recall that a slightly less standard form of the Riemann mapping theorem 
states that any doubly connected domain $\Omega$ can be conformally mapped onto an annulus of the form $\bbA_R$, where the 
number $R$ is uniquely determined by $\Omega$.  Assume that  the outer boundary component of  $\Omega$ consists of the unit 
circle and  the inner boundary component of a quite arbitrary, smooth closed curve contained in $\bbD$. 
Denoting such a conformal map by $\psi: \Omega \to \bbA_R$, its inverse 
can be taken for $\varphi$. Obviously, the class of such mappings $\varphi$  is infinitely richer than the pure M\"obius maps. 
In addition, we will assume in the next that $\varphi$ has an extension as a quasiconformal mapping on $\bbD$. 

For the formulation of the result, let us define some parameters. We fix $\beta_\infty > 0$ such that 
\bea
\max \{ |f(z)|, |f'(z)| \}  \leq \beta_\infty  \Vert f \Vert_{p,\omega}
\label{4.1}
\eea
for all $|z| \leq 1/2$,   $f \in A_\omega^p$. As well-known, the existence of such a constant follows from 
the Cauchy integral formula.
Given a quasiconformal $\varphi:\bbD \to \bbD $ with inverse $\psi = \varphi^{-1}$, we denote  
\bea 
\beta_\varphi := \max\big\{ 1, {\rm ess}\!\!\!\sup_{z \in \bbD_{1/2}} \!\!\! |\nabla \varphi(z)|, {\rm ess}\!\!\!\sup_{z \in 
\bbD_{1/2}} \!\!\! |\nabla \psi(z)| \big\} . 
\label{4.0} 
\eea
We will assume that the operators $C_\varphi$ and $C_\psi$ are bounded $L_\nu^p \to L_\nu^p$ and denote the operator norms by 
$d_\varphi$ and $d_\psi$. Recall that  $d_M = d_M(p,\omega) $ stands for the norm of the operator $M$ of Proposition 
\ref{prop 2.4} and  $d_P = d_P(p,\omega)$ for the norm of the Bergman projection $P_\omega : L_\nu^p \to A_\nu^p$, see Lemma 
\ref{lem2.0}. Finally, we let 
\bea
\delta = \min\Big\{ \frac{1}{2\beta_\varphi} ,\frac{1}{ \pi^{1/2} \beta_\varphi^{1+ p/2} 
\big( \beta_\infty d_Pd_M \max( d_\varphi,d_\psi)  \big)^{p/2} }\Big\}   \leq  \frac12 .
\label{J30}
\eea

The  result reads as follows.

\BEP  \label{th4.1} 
Let $p\in (1,2]$ and let $\omega$ be standard weight or belong to $\WA$. Let  $\varphi $ be a quasiconformal mapping from $\bbD$ 
onto itself with $\varphi(0) = 0$  such that the composition operators $C_\varphi$ and $C_\psi$, $\psi = \varphi^{-1}$,  are bounded $L_\omega^p \to L_\omega^p$ and such that $\beta_\varphi < \infty$. 
If the mapping $\varphi$ is  conformal on $\bbA_R$ for some $R \in (0, \delta )$, 
then operator  $K_\varphi= K_{\varphi,\omega} :A_\nu^p \to A_\nu^p$ is invertible.
\ENP

Proof. We denote  $S= \inf\{|\varphi(z)| \, : \, |z| \leq R \}$ and note that due to the assumption $\varphi(0)=0$,
\ef{4.0} and the mean value theorem, $S \leq \beta_\varphi R < \beta_\varphi \delta \leq 1/2$. 

Given $f \in L_\nu^p$, we follow \ef{3.14} and  write 
\bea
& & f -  K_\varphi K_\psi f =  K_\varphi \big( C_\psi - K_\psi) f
\label{4.4}
\eea
and
\bea
& & f - K_\psi K_\varphi f =  K_\psi \big( C_\varphi - K_\varphi \big) f   \label{4.6}
\eea
and next show that 
the norms of these two expressions are at most 
$C \Vert f \Vert_{p,\omega}$ for a constant  $C>0$ smaller than one, which implies the claimed invertibility of $K_\varphi$ 
in the same way as at the end of the proof of Theorem \ref{th3.1}. 

In order to use Proposition \ref{prop 2.4}, let us next show that $f \circ \psi \in H$. Since the composition
operator is assumed to be bounded in $L_\nu^p$, there holds $f \circ \psi \in L_\nu^p$. In addition, we have
$\bar \partial (f \circ \psi) = ( \partial f \circ \psi) \bar \partial \psi$ and 
\bea
& &  \Vert \bar \partial ( f \circ \psi)  \Vert_{p,\omega}^p
= \int\limits_{\varphi (\bbA_R)} | \partial f \circ \psi |^p \, | \bar \partial \psi|^p \omega dA +
\int\limits_{\varphi (\bbD_R)} | \partial f \circ \psi |^p \, | \bar \partial \psi|^p \omega dA
\roweq 
0 + 
\int\limits_{\varphi(\bbD_R)} | \partial f \circ \psi |^p \, | \bar \partial \psi|^p \omega dA.
\label{4.12}
\eea
By $R,S < 1/2$, \ef{4.1}, \ef{4.0}, this expression is bounded by 
\bea
\max \limits_{z \in \bbD_{1/2}} |f'(z)|^p \int\limits_{\bbD_{S}}  | \nabla \psi|^p \omega dA
\leq \pi S^2  \beta_\infty^p \beta_\varphi^p  \Vert f \Vert_{p,\omega}^p , \label{4.12g}
\eea
This implies that $\bar \partial (f \circ \psi) \in L_\nu^p$ and thus $f \circ \psi \in H$. Hence, 
Proposition \ref{prop 2.4} and \ef{4.12}--\ef{4.12g} yield the estimate 
\bea
& & \Vert C_\psi f  - K_\psi f\Vert_{p, \omega}^p 
= \Vert (P_\omega - I)C_\psi  f  \Vert_{p,\omega}^p 
= \Vert Q_\omega  ( f \circ \psi)  \Vert_{p,\omega}^p
\roweq 
\Vert M \bar \partial ( f \circ \psi)  \Vert_{p,\omega}^p
\leq  
d_M^p \Vert \bar \partial ( f \circ \psi)  \Vert_{p,\omega}^p
\leq \pi S^2  \beta_\infty^p \beta_\varphi^p d_M^p \Vert f \Vert_{p,\omega}^p. 
\label{4.12gg}
\eea
We obtain   
\bea
& &  \Vert K_\varphi \big( C_\psi   - K_\psi \big) f\Vert_{p, \omega}
\leq  (\pi S^2)^{1/p} d_P d_\varphi d_M  \beta_\infty \beta_\varphi  \Vert f \Vert_{p,\omega}^p. 
\label{4.12ggg}
\eea
Since  $S < \beta_\varphi \delta$, the coefficient on the right is strictly smaller than
\bea
\delta^{2/p} \pi^{1/p} \beta_\varphi^{1+ 2/p} \beta_\infty d_P d_\varphi d_M    \leq 1, 
\label{4.14}
\eea 
due to the choice of $\delta$ in \ef{J30}. Thus, from \ef{4.4}, \ef{4.12ggg}--\ef{4.14} we get that $K_\varphi K_\psi$ is invertible.

The expression \ef{4.6} can be estimated similarly to \ef{4.12}--\ef{4.12ggg} with obvious changes, e.g., replacing the
integration domains in \ef{4.12} by $\bbA_R$ and $\bbD_R$. This   yields the bound
\beas
 \Vert K_\psi \big( C_\varphi f  - K_\varphi \big)  f\Vert_{p, \omega}
\leq (\pi R^2)^{1/p} d_Pd_M d_\psi\beta_\infty \beta_\varphi  \Vert f \Vert_{p,\omega}^p 
=: C \Vert f \Vert_{p,\omega}, 
\eeas
and as above, we deduce that also $K_\psi K_\varphi$ is invertible. 
This completes the proof. \ \ $\Box$

\section{Examples} \label{sec5}

In the next examples we consider projected composition operators in the simplest case of the standard 
unweighted Bergman-Hilbert space $A^2$.

\BEE
Let us define the mapping $\varphi$ on $\bbD$ as
\bea
\varphi(z) = \varphi(re^{i \theta}) = re^{i\theta + i b(r) }= z e^{ i b(r) }   \label{J0}
\eea
where $b: [0,1) \to \bbR $ is a continuously differentiable function. We denote $b'(r) = db(r) /dr$ and claim that if
\bea
|b'(r)| \leq C ( 1- r^2), \ \ \ \mbox{where} \ C = 
\min(1, \gamma_\psi, \gamma_\varphi).  \label{J1}
\eea
then $\varphi$ is quasiconformal and  satisfies the assumptions of Theorem \ref{th3.1}
and $T_\varphi:A^2 \to A^2$ is thus invertible. Here, $\gamma_\psi$ and $\gamma_\varphi$ are as in
\ef{3.2}, \ef{3.2a}.

We calculate
\bea
& & \partial \varphi(z) = \partial(z e^{i b(r)}) = e^{i b(r) } + iz  e^{i b(r) } \partial b(r) = 
 e^{i b(r) }\Big(1 +  \frac{i}2 r b'(r) \Big) , 
\row
\bar \partial \varphi(z) = z \bar \partial e^{i b(r)} = 
\frac{i z^2}{2 r }  b'(r)  e^{i b(r) }
\eea
This yields 
\bea
& & |\mu(z)| = \frac{|\bar \partial \varphi(z)|}{|\partial \varphi(z)|} =   
\frac{r |b'(r)| }{ 2|1 + \frac{1}2 r b'(r) |}  \leq \frac{C r (1-r^2)}{ 2|1 - C r(1-r^2) /2 |}.  \label{J2}
\eea
Here, we proceed with the bounds $r(1-r^2) \leq 2/(3 \sqrt{3)} \approx 0.3849$ and $C \leq 1 $, which show that \ef{J2} is not larger than
\bea
& & \frac{ Cr(1-r^2)}{2(1- 1/(3\sqrt{3}) ) }  \leq 
\frac{C}2 (1-r^2) \frac{3\sqrt{3}}{3\sqrt{3} -1 }   
\row<
\frac12 \min (1, \gamma_\psi, \gamma_\varphi) (1-r^2) \frac{3\sqrt{3}}{2\sqrt{3} }  
= 
\frac34  \min (1, \gamma_\psi, \gamma_\varphi) (1-r^2). \label{J4}
\eea
This shows that $\varphi$ is quasiconformal. 
Now, there holds $|\varphi(z)| = |z|$ and, by the choice of the space,  $\alpha =0$ so that the right-hand side of \ef{3.2} 
reduces to $\gamma_\psi (1-|z|^2)$. Thus, \ef{J2}--\ef{J4} imply that \ef{3.2} holds. Also, since $|\psi (z)| = |\varphi^{-1} (z)| = 
|z e^{- i b(r)}| = |z|$, condition \ef{3.2a} is satisfied as well, due to \ef{J4}.  We conclude that \ef{J1}
is sufficient for the invertibility of $T_\varphi$. 
\ENE

\BEE
Let $a \in (0, \infty)$ and $R\in (0,1)$ and let us define 
\begin{equation}
\varphi(z) = \left\{
\begin{array}{ll}
R^{1-a} z |z|^{a-1} , \ \ &\mbox{if} \ |z| \leq  R,  \\
z &\mbox{for $R < |z| < 1$} . 
\end{array}
\right.  \label{J6}
\end{equation}
Then, $\varphi$ is weakly differentiable and belongs to the Sobolev space $W^{1,2}$.
According to Example 10.1.(2) in \cite{Kos}, $\varphi$ is quasiconformal, since its 
Beltrami coefficient satisfies $|\mu(z)| = (a-1)/(a+1)$ for $|z| \leq R$ and, obviously, $\mu(z) = 0$ for $|z| > R$.
By a direct calculation, 
\begin{equation}
\psi(z) := \varphi^{-1}(z) = \left\{
\begin{array}{ll}
R^{1-1/a} z |z|^{1/a-1} , \ \ &\mbox{if} \ |z| \leq  R,  \\
z &\mbox{for $R < |z| < 1$} . 
\end{array}
\right.
\end{equation}
We claim that if 
\bea
a < \frac{1+ \min(\gamma_\psi,\gamma_\varphi)(1-R^2)}{1- \min(\gamma_\psi,\gamma_\varphi)(1-R^2)}  \label{J10}
\eea
then  $T_\varphi$ is invertible. Since $\varphi$ is analytic in $\{|z| > R\}$, it suffices to prove \ef{3.2}, \ef{3.2a} 
for $|z| \leq R$. Now, \ef{J10} implies
\bea
|\mu(z)| = \frac{a-1}{a+1} 
< \min(\gamma_\psi,\gamma_\varphi) (1-R^2)  , \label{J12}
\eea
and, in view of the facts $|\varphi(z)|^2 = R^{2-2a}|z|^{2a}$ and $|\psi(z)|^2 = R^{2-2/a}|z|^{2/a}$, the expression on the right is not larger than
\bea
\gamma_\psi  (1-|\varphi(z)|^2) \ \ \ \mbox{or} \ \ \ 
\gamma_\varphi (1-|\psi(z)|^2)  
\eea
for $|z| \leq R$. We thus obtain \ef{3.2}, \ef{3.2a} from \ef{J12}.
\ENE

\BEE
Let us complement the above considerations by presenting an example of a non-invertible projected composition 
operator in $A^2$ such that the symbol is quasiconformal on $\bbD$ and 
equal to  the identity on an  annulus $\bbA_R$. For simplicity, we start by introducing a symbol which is not 
quasiconformal, namely, we choose $R = (3/7)^{1/4} \approx 0.8091$ and define
\begin{equation}
{\widetilde \varphi}(z) = \left\{
\begin{array}{ll}
z , \ \ &\mbox{if} \ |z| > R,  \\
- Re^{i \theta} &\mbox{for $ z = re^{i \theta}$ with $r \leq R$} . 
\end{array}
\right.
\end{equation}
This can also be written as 
\bea
{\widetilde \varphi}(r e^{i \theta}) = \tilde b(r) e^{i \theta + i\tilde a(r)} ,   \ \ \ r e^{i \theta} \in \bbD,
\label{4.26}
\eea
where
\begin{equation}
\tilde a(r)  = \left\{
\begin{array}{ll}
\pi , \ \ &\mbox{if} \ r\in [0,R] ,  \\
0 &\mbox{if} \ r\in [R , 1)]    
\end{array}
\right. \ \ \ \mbox{and}  \ \ \ 
\tilde b(r)  = \left\{
\begin{array}{ll}
R , \ \ &\mbox{if} \ r\in [0,R] ,  \\
r &\mbox{if} \ r\in [R , 1)]    
\end{array}
\right.  \label{4.28}
\end{equation} 
Apparently, the restriction of ${\widetilde \varphi}$ onto $\bbA_R$ is conformal. However, the operator $K_{\widetilde \varphi}$ is
not invertible in $A^2$, since it maps the function $f(z) = z $ into 0. To see this, we calculate
\bea
& & P C_{\widetilde \varphi} f (z)= \frac1\pi \int\limits_0^{1} \int\limits_0^{2 \pi} \frac{\tilde b(r)  
e^{i\theta + i\tilde a(r)}}{(1- z re^{-i\theta})^2 }  rd\theta dr 
\roweq
\frac1\pi \int\limits_0^{1} \int\limits_0^{2 \pi} \sum_{n=0}^\infty (n+1) z^n \tilde b(r) r^{n+1} e^{i \theta (1-n)} d\theta
e^{i a(r)} dr = 4 z \int\limits_0^{1}  \tilde b(r)r^2 e^{i a(r)} dr  . \label{4.29}
\eea

The expression \ef{4.29} vanishes, since it is $4z$ times 
\bea
& & \int\limits_0^{R}  \tilde b(r)r^2 e^{i \tilde a(r)} dr +
 \int\limits_R^{1}  \tilde b(r)r^2  e^{i \tilde a(r)} dr
\roweq
 - R \int\limits_0^{R} r^2  dr  +  2 \int\limits_{R}^1 r^3 dr'= 2 \Big( -
\frac{R^4}3 + \frac{1 - R^4}{4}	\Big) = 0  \label{4.30}
\eea
due to the choice of $R$. 

The rest of the example consists of a technical modification of the above symbol ${\widetilde \varphi}$ so 
as to make it quasiconformal while still preserving the non-injectivity of the corresponding projected 
composition operator. This can be achieved by defining (cf. \ef{4.26}--\ef{4.28})
\bea
\varphi(r e^{i \theta}) = b(r) e^{i \theta + i a(r)} ,   \ \ \ r e^{i \theta} \in \bbD,
\eea
where the smooth functions $a: [0,1) \to [0, \pi]$ and $b:[0,1)\to [0,1]$ (modifications of 
$\tilde a$ and $\tilde b$) will be  chosen such that 
$\varphi$ becomes a diffeomorphism of $\bbD$ onto itself. Then, the quasiconformality of $\varphi$ 
follows from \cite{Kos}, Example 1.1.5., since $\varphi$ is analytic in a neighborhood of $\partial \bbD$. 
The details will require some efforts. 

To define $a$, we let $R =(3/7)^{1/4} \approx 0.8091$ be as above and set $R' = 9/10$. Then, we consider arbitrary parameters 
$\delta_a \in (0,1/20)$  and $\delta \in (0,1/20)$ and set 
\begin{equation}
a(r)  = \left\{
\begin{array}{ll}
\pi , \ \ &\mbox{if} \ r\in [0,R] ,  \\
a_1(r) 
, \ \ &\mbox{if} \ r\in [R,R+ \delta] ,  \\
0  , \ \ &\mbox{if} \ r\in [R + \delta ,R'] ,  \\
a_2 (r) , \ \ &\mbox{if} \ r\in [R' ,R' + \delta_a] ,  \\
0 &\mbox{if} \ r\in [R' + \delta_a, 1)]  .    \label{4.34}
\end{array}
\right.
\end{equation}
To choose the function $a_1$, we first take a $C^\infty$-function $\widehat a_1: \bbR \to [0,\pi]$ such that 
$\widehat a_1(r) = \pi$
for $r \leq 0$ and $\widehat a_1(r) = 0$ for $r \geq 1$ and such that it decreases
monotonically on the subinterval $ (0,1)$. We set $a_1(r) = \widehat{a}_1( (r - R)/\delta))$. Then, 
$a_1$ becomes a  monotonically decreasing  $C^\infty$-function with  $a_1(R) = \pi$ and $a_1(R+ \delta)=0$.

As for the function $a_2$, we first choose a $C^\infty$-function $\widehat a_2: \bbR \to [-1,0]$ with $\widehat a_2(-1)=
\widehat a_2(1) = 0$, $\widehat a_2(0) = -1$ and  $\widehat a_2(r)= 0$  for $|r| \geq 1$ such that it  decreases
monotonically on the subinterval $ (-1,0)$ and increases
monotonically on the subinterval $ (0,1)$. Then, we define $a_2$ by scaling: we set
$$
a_2(r) = \widehat{a}_2\Big( \frac{2(r - R')}{\delta_a} -1 \Big ) . 
$$ 
Consequently, $a_2: [R',R'+\delta] \to (-\infty, 0]$ is a $C^\infty$-function  decreasing monotonically 
from the value $\widehat{a}(R') =  0 $ to the
value $\widehat{a}( R' + \delta_a/2) = -1$ and then increasing monotonically to the value $\widehat{a}(R' + \delta_a)=0$,
and the values of $a_2$ are  negative on $(R',  R' +  \delta_a)$. It is also clear that with these definitions, the
function $a$ in \ef{4.34}  becomes $C^\infty$-smooth.

We choose $b$ such that 
\begin{equation}
b(r)  = \left\{  
\begin{array}{ll}
b_1(r)  , \ \ &\mbox{if} \ r\in [0, R] ,  \\
r, \ \ &\mbox{if} \ r \in [R,1]     .  \label{4.36}
\end{array}
\right.
\end{equation}
where $b_1 : [0,R] \to [0,1]$ is a smooth, increasing function such that $b_1(0) = 0$ and 
$b_1(R) =R$; we will fix $b_1$ later.

In the same way as in \ef{4.29} we see that $PC_\varphi f(z) = 4 z \int\limits_0^{1}  b(r)r e^{i a(r)} dr$,
where 
\beas
& & \int\limits_0^{1}  b(r)r^2 e^{i a(r)} dr = \int\limits_0^{1}  b(r)r^2 \cos( a(r)) dr
+ i \int\limits_0^{1}  b(r)r^2 \sin( a(r)) dr 
\roweq
\int\limits_0^{1}  b(r)r^2 \cos( a(r)) dr + i \int\limits_R^{R'+ \delta_a}  r^3 \sin( a(r)) dr = :I_{\rm Re}(\delta_a, \delta) +
i I_{\rm Im} (\delta_a, \delta) ,
\eeas
due to the choice of the functions $a$ and $b$. 

We now note  that for every sufficiently small $\epsilon > 0$, we can choose 
$\delta_a \in (0, \epsilon)$ and $ \delta \in (0, \epsilon)$ such that $I_{\rm Im} (\delta_a, \delta)= 0$ $(*)$. This follows from the observation that 
\bea
& & I_{\rm Im} (\delta_a, \delta) = \int\limits_R^{R+ \delta}  r^3 \sin( a(r)) dr
+\int\limits_{R'}^{R'+ \delta_a}  r^3 \sin( a(r)) dr =: J_1 (\delta) + J_2 (\delta_a),   \label{4.37}
\eea
where, due to the choice of the function $a_1$ by scaling with respect to $\delta$, the function 
$\delta \mapsto J_1(\delta)$ is continuous and positive on the interval   $(0,1/20)$ and 
satisfies $\lim_{\delta \to 0} J_1(\delta) = 0$. Similarly, the function $\delta_a \mapsto J_2(\delta_a)$ is continuous and 
negative on the interval  $(0,1/20)$ and has  the property  $\lim_{\delta_a \to 0} J_2(\delta_a) = 0$. This implies 
the claim about the choice $(*)$.

In order to choose $b_1$, we consider another parameter $\delta_b \in (0,1/10)$ and, first, a piecewise linear, monotonely 
increasing function $\widehat b$ such that  $\widehat b (r) = r$ for $0 \leq r < \delta_b/2$, $\widehat b (\delta_b) = R - 
\delta_b$ and $\widehat b (r ) = r$ for $r \geq R$.  It is not difficult to see that, using 
a suitable $C^\infty$  mollifier (e.g. $\Phi_{\delta_b/10}$, see the remark after the proof) and defining $b_1$ as the convolution  of the mollifier and $\widehat b$, 
the function
\beas
\delta_b \mapsto  - \int\limits_0^R b_1(r) r^2  dr =: K_1(\delta_b)
\eeas
is continuous on the interval $(0,1/10)$ and that 
\beas  
\lim_{\delta_b \to 0} K_1(\delta_b) = - \int\limits_0^R R  r^2  dr = B_1.
\eeas
(Note also that $b : [0,1) \to [0,1)$ becomes an invertible $C^\infty$-function.)
We also have $K_1(\delta_b) > B_1$ for all $\delta_b \in (0,1/10)$. In particular, there exists $\epsilon_0>0$ such that for all $\epsilon' \in (0,\epsilon_0) $ we can find $\delta_b$ such that
\bea
K_1(\delta_b) -  B_1 = \epsilon' . \label{4.38}
\eea

As the last step we note that, for any choice of the function $a(r)$ as above, we have 
\beas
K_2(\delta_a, \delta) := \int\limits_R^1 b(r) r^2 \cos (a(r)) dr = \int\limits_R^1 r^3 \cos (a(r)) dr
< \int\limits_R^1  r^3 dr = B_2.  
\eeas
and, due to the definition of $a$, the difference 
\beas
B_2 - K_2(\delta_a, \delta)   
\eeas
can be made arbitrarily small, if $\delta_a$ and $ \delta$ are small enough. 
Accordingly, we finally fix $\delta_a, \delta$ such that  $I_{\rm Im} (\delta_a, \delta) = 0$ (see $(*)$) 
and such that  $B_2 - K_2(\delta, \delta_a)  < \epsilon_0$. 
We then choose $\delta_b$ such that \ef{4.38} holds with $\epsilon' = B_2 - K_2(\delta_a, \delta)$. This yields
\beas
& & I_{\rm Re} (\delta_a, \delta) = \int\limits_0^1 b(r) r^2 \cos(a(r)) dr 
= K_1(\delta_b)  + K_2(\delta_a, \delta)  =  B_1 + \epsilon' + B_2 - \epsilon' = 0 ,
\eeas
by  \ef{4.30}. This shows that $PC_\varphi f = 0$ and thus $K_\varphi$ is not invertible. 

We finally remark that $\varphi$ is a diffeomorphism and thus quasiconformal, since its inverse 
is the smooth function $re^{i \theta}  \mapsto b^{-1}(r) e^{i \theta - i a(b^{-1} (r))} $. Note that $\varphi(z) =z$
in a neighborhood of 0, by the choice of the functions $a$ and $b$. \ \ $\Box$
\ENE

\bigskip

The standard mollifier $\Phi_\alpha$, $0 < \alpha \leq 1$, is defined as $\Phi_\alpha (x) = \Phi(x/\alpha)/\alpha$, where
$\Phi(x) = C\exp(-1/(1-x^2))$ for $|x| < 1$ and $\Phi(x) = 0 $ for $|x| \geq1$ and the normalization factor $C > 0$
makes the integral of $J$ equal one.

\subsection*{Acknowledgment} The first named author was supported  by a T\"UBITAK project with project number 
1059B192300654. The second named author acknowledges the support of the Academy of Finland projects no. 359563
and 359642.


\begin{thebibliography}{1}

\bibitem{A} L. Ahlfors, Lectures on quasiconformal mappings, Mathematical Studies vol.10, D. Van Nostrand Company 
Inc., Princeton NJ (1966). 

\bibitem{AS} A. Aleman, A. Siskakis, Integration operators on Bergman spaces,  Indiana Univ. Math. J. 46, 2 (1997), 337--356.

\bibitem{APP} H. Arroussi, I. Park, J. Pau, Schatten class Toeplitz operators acting on large weighted
Bergman spaces, Studia Math. 229, 3 (2015), 203--221.

\bibitem{ATV} H. Arroussi, J. Taskinen, J. Virtanen, Estimates for the derivatives of the Bergman kernel in large Bergman spaces and 
applications to operator theory, submitted. 

\bibitem{Ast} K. Astala, A remark on quasiconformal mappings and BMO-functions. Michigan Math. J. 30, 209–12 (1983) 

\bibitem{BB} B. Berndtsson, Weighted estimates for $\bar \partial$ in domains in $\bfC$.
Duke Math. J. 66 (1992), no. 2, 239--255.


\bibitem{const2010} O. Constantin,
Carleson embeddings and some classes of operators on weighted Bergman spaces.
J. Math. Anal. Appl. 365 (2010), no. 2, 668-682. 
	
\bibitem{Cuckovic} Z. Cuckovic,  Projected composition operators on pseudoconvex domains, Integr. Eq. Oper. Theory 93 (35), (2021).	

\bibitem{D2} M. Dostani\'c, {Integration operators on Bergman spaces with exponential weights}, Revista Mat. Iberoamericana 23 (2007), 421--436. 

\bibitem{Fang} X. Fang, K. Guo, Z. Wang, Fang, Composition operators on the Bergman space with quasiconformal symbols,
J. Geometric Anal. (2023) 33:125

\bibitem{FS}  J. E. Fornaess and N. Sibony, On $L\sp p$ estimates for $\overline\partial$, Several complex variables and complex geometry, Part 3 (Santa Cruz, CA, 1989), Proc. Sympos. Pure Math., vol. 52, Amer. Math. Soc., Providence, RI, 1991, pp. 129--163.
	

\bibitem {GP1}  P. Galanopoulos, J. Pau, {Hankel operators on large weighted Bergman spaces}, Ann. Acad. Sci. Fenn. Math. 37 (2012), 635--648.



\bibitem{GS} N.G. G\"oğ\"uş, N., S.Y. S\"onmez, Toeplitz Operators on Weighted Bergman Spaces on Finitely Connected Domains. 
Annals of Functional Analysis (2022)	


\bibitem{Ho}  L. H\"ormander, An introduction to complex analysis in several variables, North-Holland Mathematical Library, vol. 7, North-Holland Publishing Co., Amsterdam, 1990.

\bibitem{HoN}  L. H\"ormander, Notions of convexity BOOK. 


\bibitem{HLS} Hu, Z., Lv, X., Schuster, A.: Bergman spaces with exponential weights. J.Functional Anal. \textbf{276}, (2019), 1402–1429


\bibitem{HP2021} Z. Hu and J. Pau, Hankel operators on exponential Bergman spaces. Sci. China Math. 65, 421--442 (2022). 

\bibitem{Kos} P. Koskela, Lectures on quasiconformal and quasisymmetric mappings, University of Jyv\"askyl\"a (2009).

\bibitem{Kos14} P. Koskela, J. Xiao, Y. Ru-Ya Zhang, Y.  Zhou,  
A quasiconformal composition problem for the Q-spaces, J. Eur. Math. Soc. 19, 1159--1187


\bibitem{PP1}  J. Pau, J.A. Pel\'{a}ez, {Embedding theorems and integration operators on Bergman spaces with rapidly decreasing weights}, J. Funct. Anal. 259 (2010), 2727--2756.

\bibitem{PP2} J. Pau, J.A. Pel\'{a}ez, {Volterra type operators on Bergman spaces with exponential weights}, Contemp. Math. 561 (2012), 239--252.

	
\bibitem{Mishra} S. Pattanayak, C.K. Mohapatra, A.K. Mishra, A new class of composition operators, Internat J. Math \& Math. Sci. \textbf{10} (1987), 473-482.

\bibitem{PRAdv} J.A. Pel\'aez, J. R\"atty\"a, Bergman projection induced by radial weight, 
Advances Math. 391(2021) 107950 1--70.

\bibitem{Range} M. Range, Holomorphic functions and integral representations in several complex variables,
 Graduate texts in mathematics 108 Springer New York Berlin Heidelberg (1986)
 
\bibitem{Rei} H.M. Reimann, Functions of bounded mean oscillation and quasiconformal mappings. Comment. Math. Helv. 49 (1974), 260–276 


\bibitem{Rochberg} R. Rochberg, Projected Composition Operators on the Hardy Space,  Indiana Univ. Math. J. 43 (2), (1994).

\bibitem{ST} S. S\"onmez, J. Taskinen, Relaxing the positivity assumption on the symbol of a
Bergman-Toeplitz operator. To appear in Proc. Amer. Math. Soc.


\bibitem{Vai} J. Väisälä, Lectures on n-dimensional quasiconformal mappings. Lecture Notes in Mathematics, 
Springer-Verlag, Berlin, 1971.

\bibitem{Zhao} C. Zhao, Boundedness of projected composition operators over the unit disc,  J. Math. Anal. Appl. 467 (2018), 521-536.

\bibitem{Zhu}  K. Zhu, Operator Theory in Function Spaces, 2nd ed. Vol. 138, 
Mathematical surveys and monographs, American Mathematical Society, Providence 
RI, 2007.




\end{thebibliography}
\end{document}